\documentclass[11pt,twoside]{article}
\usepackage{latexsym}
\usepackage{amssymb,amsbsy,amsmath,amsfonts,amssymb,amscd}
\usepackage{graphicx,dsfont}
\usepackage{hyperref}
\usepackage{caption,color}
\newcommand{\commentout}[1]{}
\newcommand{\R}{\mathbb{R}}

\newcommand {\fer}   {\eqref}

\newcommand {\al} {\alpha}
\newcommand {\e}  {\varepsilon}

\newcommand {\vp} {\varphi}

\newcommand {\Chi} {{\bf \raise 2pt \hbox{$\chi$}} }
\newcommand {\caa} { {\mathcal A} }

\newcommand {\f}   {\frac}
\newcommand {\p}   {\partial}

\newcommand {\proof} {\noindent {\bf Proof}. }
\newcommand{\beq}{\begin{equation}}
\newcommand{\eeq}{\end{equation}}
\newcommand{\bal}{\begin{align}}
\newcommand{\bc}{\begin{cases}}
\newcommand{\ec}{\end{cases}}
\newcommand{\bea} {\begin{array}{rl}}
\newcommand{\eea} {\end{array}}
\newcommand{\bepa}{\left\{ \begin{array}{l}}
\newcommand{\eepa} {\end{array}\right.}
\newtheorem{theorem}{Theorem}[section]
\newtheorem{lemma}[theorem]{Lemma}

\newtheorem{definition}[theorem]{Definition}

\newtheorem{prop}[theorem]{Proposition}

\newcommand{\qed}{{ \hfill
                     {\unskip\kern 6pt\penalty 500 \raise -2pt\hbox{\vrule\vbox to 6pt{\hrule width 6pt
                     \vfill\hrule}\vrule} \par}   }}
\title{\Large \bf Singular limits for models of selection and mutations with heavy-tailed mutation distribution }

\author{
  Sepideh Mirrahimi\thanks{Institut de Math\'ematiques de Toulouse; UMR 5219, Universit\'e de Toulouse; CNRS, UPS, F-31062 Toulouse Cedex 9, France; E-mail: \texttt{Sepideh.Mirrahimi@math.univ-toulouse.fr}}}

\date{\today}

\begin{document}
\maketitle
\pagestyle{plain}
\pagenumbering{arabic}

\begin{abstract}
In this article, we perform an asymptotic analysis of a nonlocal reaction-diffusion equation, with a fractional laplacian as the diffusion term and with a nonlocal reaction term. Such equation models the evolutionary dynamics of a phenotypically structured population. \\
We perform a rescaling considering large time and small effect of mutations, but still with algebraic law.  With such rescaling, we expect that the phenotypic density will concentrate as a Dirac mass which evolves in time. To study such concentration phenomenon, we extend an approach based on Hamilton-Jacobi equations with constraint, that has been developed to study models from evolutionary biology,  to the case of fat-tailed mutation kernels. However, unlike previous works within this approach, the WKB transformation of the solution does not converge to a viscosity solution of a Hamilton-Jacobi equation but to a viscosity supersolution of such equation which is minimal in a certain class of supersolutions.  Such property allows to derive the concentration of the population density as an evolving Dirac mass, under monotony conditions on the growth rate, similarly to the case with thin-tailed mutation kernels.

\end{abstract}

\noindent\textbf{AMS Class. No:}  35K57, 35B25, 47G20, 49L25,  92D15 

\section{Introduction}
\label{sec:intro}

\subsection{Model and motivation}

In this paper we are interested in the following selection-mutation model
\beq
\label{SM}
\begin{cases}
\p_t n  +  (-\Delta)^{\al} n = n\, R(x,I),\\
n(x,0)=n^0(x),\ x\in {\mathbb{R}^d},
\end{cases}
\eeq
with 
\beq
\label{def:I}
I(t)=\int_{{\mathbb{R}^d}}  n(t,x)dx.
\eeq
In all what follows,  $\alpha\in (0,1)$ is given. The term $(-\Delta)^{\al}$ denotes the fractional {laplacian}:
\beq
\label{FL}
(-\Delta)^{\al} n(t,x) = \int_{h\in \R^d} \left[ n(t,x) - n(t,x+h)  \right] \f{dh}{|h|^{d+2\al}}.
\eeq
Equation \fer{SM} has been derived from a stochastic individual based model describing the evolutionary dynamics of a phenotypically structured population \cite{BJ.SM.WW:11}. Here, $t$ corresponds to time and $x$ corresponds to a phenotypic trait. The function $n$ represents the {phenotypic density} of a population. The term $I(t)$ corresponds to the total population size. The growth rate of the individuals is denoted by $R(x,I)$ which depends on the phenotypic trait and the total population size, taking into account in this way  competition between the individuals. The fractional laplacian term models the mutations. The choice of a fractional laplacian rather than a classical laplacian or an integral kernel with thin tails, allows to take into account large mutation jumps with a high rate \cite{BJ.SM.WW:11}.
\\

\noindent
Several frameworks have been used to study models from evolutionary biology. Game theory is one of the first approaches which has contributed a lot to the understanding of mechanisms of evolution \cite{JM:74,JH.KS:88}. Adaptive dynamics, a theory  based on stability analysis of dynamical systems, allows to study evolution under very rare mutations \cite{SG.EK.GM.JM:97,OD:04}. Integro-differential models are used to study evolutionary dynamics of large populations (see for instance \cite{PM.GW:00,AC.SC:04,OD.PJ.SM.BP:05,LD.PJ.SM.GR:08}). Probabilistic tools allow to study populations of small size \cite{NC.AM:07} and also to derive the above models in the limit of large populations \cite{NC.RF.SM:08}.\\

\noindent
Within the integro-differential framework, an approach based on Hamilton-Jacobi equations with constraint has been developed during the last decade to study asymptotically, in the limit of small mutations and large time, integro-differential models from evolutionary biology. There is a large literature on this approach which was first suggested by \cite{OD.PJ.SM.BP:05}. See for instance \cite{GB.BP:08,GB.SM.BP:09,AL.SM.BP:10} where the basis of this approach for models from evolutionary biology were established. Note  that this approach has also been used  to study the propagation phenomena in local reaction-diffusion equations (see for instance \cite{MR833742,MF:85,CE.PS:89,GB.LE.PS:90}). 
The present article follows an earlier work \cite{SM.SM:15} which was an attempt to extend the Hamilton-Jacobi approach to the case where the diffusion is modeled by a fractional laplacian rather than a classical laplacian  or an integral kernel with thin tails.\\

\noindent
We consider a rescaling introduced in  \cite{SM.SM:15}, rescaling the size of the mutations to be smaller and performing a change of variable in time, to be able to observe the effect of small mutations on the dynamics. To this end, we choose $k>0$ and $\nu\in S^{d-1}$ such that $h=(e^k-1)\nu$, and perform the following rescaling
\begin{equation}
\label{res-2}
  t\mapsto \f{t}{\e}, \quad    \widetilde{M}(h,dh)=\f{dh}{|h|^{d+2\al}}=\f{e^{k } dkdS}{|e^{k }-1|^{1+2\al}} \mapsto M_\e(k,dk,dS)= \f{e^{\f k \e} \f{dk }{\e}dS}{|e^{\f k\e }-1|^{1+2\al}}.
\end{equation}
{Note that with the above transformation, the dimension $d$ disappears in the power of $|e^k -1|$, since $dh$ becomes $(e^k-1)^{d-1}e^kdkdS$. }\\
{We then study the following rescaled problem:}
 \beq
\label{SMe}
\begin{cases}
\e \p_t n_\e (t,x) = \int_0^\infty \int_{\nu\in S^{d-1}} \left( n_\e(t,x+(e^{\e k}-1)\nu) - n_\e(t,x) \right) \f{e^k }{|e^k-1|^{1+2\al}}dSdk+n_\e(t,x) \,R(x,I_\e(t)),\\
I_\e(t)=\int_{\mathbb{R}^d}  n_\e(t,x)dx,\\
n_e(x,0)=n_\e^0(x).
\end{cases}
\eeq
{With this rescaling, we consider much smaller mutation steps. The mutations' distribution has still algebraic tails but with a large power.} In particular, it has a finite variance of order $\e^2$. {Note indeed that the covariance matrix $\mathbf{v}=(v_{i,j})_{1\leq i,j\leq d}$ of the mutations' distribution above  is given by 
$$
v_{i,j}= \int_0^\infty \int_{\nu\in S^{d-1}} (e^{\e k}-1)^2  \nu_i \nu_j \f{e^k }{|e^k-1|^{1+2\al}}dSdk=O(\e^2), \quad \text{as $\e\to 0$}.
$$
}

\noindent 
The rescaling  \fer{res-2} is very different from the one considered   for a model with a classical laplacian \cite{GB.BP:08,GB.SM.BP:09}, that is 
$$
\e \p_t  n_\e -\e^2 \Delta n_\e= n_\e R(x,I_\e),
$$
or the one considered for a model with an integral kernel $J$ with thin tails \cite{GB.SM.BP:09}, that is
\beq
\label{eq:BMP}
\e \p_t  n_\e - \int_{\R^d} \big( n_\e(t,x+h)-n_\e(t,x) \big) J \big(\f{h}{\e} \big)\f{dh}{\e^d}= n_\e R(x,I_\e).
\eeq
The possibility of big jumps in \fer{SM} changes drastically the behavior of the solutions and leads to much faster dynamics of the phenotypic density. {Therefore, such type of rescaling cannot be used. In particular, if we followed the same method that has been used in \cite{GB.SM.BP:09} for \fer{eq:BMP}, to study \fer{SM} we would obtain a Hamiltonian, at the limit as $\e\to 0$,  that has infinite value (see equation (16) in  \cite{GB.SM.BP:09}).  To obtain a relevant equation at the limit and similar type of behavior as in \cite{GB.SM.BP:09} we consider a rescaling in  \fer{res-2} that makes the size of the mutations much smaller}.  The rescaling  \fer{res-2} is derived thanks to an analogy to the fractional Fisher-KPP equation \cite{SM.SM:15}. In \cite{SM.SM:15}, an asymptotic analysis was provided in the case of the fractional Fisher-KPP equation where the propagation has an exponential speed \cite{XC.AC.JR:12,XC.JR:13} leading to significantly different scalings {compared to the case of the classical Fisher-KPP equation (see for instance \cite{CE.PS:89})}. Model \fer{SMe} was then derived with an inspiration from such rescaling. Note however that in all of the above rescalings the variance of the rescaled mutation kernel is of order $\e^2$. To be able to observe concentration phenomena,  the variance of the mutation kernel must be indeed small.\\

\noindent{
An asymptotic analysis of \fer{SMe} was provided   in \cite{SM.SM:15}  for homogeneous reaction terms $R(I)$ and under strong assumptions on the initial data. Here, we extend this result to the case of heterogeneous $R(x,I)$ and relax the assumptions on the initial data, obtaining in this way a result which is analogous to the previous works with standard terms of mutation \cite{GB.BP:08,GB.SM.BP:09}.}  \\

\noindent
The method developed in \cite{SM.SM:15} has been extended in several directions. In \cite{PS.AT:18,AL:18} an asymptotic study of a Fisher-KPP type equation has been provided in periodic media and with a general non-local stable operator of order $\alpha \in (0,2)$. In \cite{EB.JG.CH.FP:18}  a homogeneous Fisher-KPP type  model has been studied,  modeling the diffusion  by a convolution term without singularity  but considering more general  decays for the integral kernel. The method provided in the present paper can also be used to generalize the results of \cite{EB.JG.CH.FP:18} and to study selection-mutation models with the integral kernels  given in \cite{EB.JG.CH.FP:18}, where a similar difficulty appears.

\subsection{Assumptions}

\noindent
Before presenting our assumptions, we first introduce the classical Hopf-Cole transformation
\beq
\label{Hopf}
n_\e=\exp\left( \f{u_\e}{\e} \right).
\eeq
Here are our assumptions:\\
We  assume that there are two constants $0<I_m<I_M<\infty$ such that
\begin{equation}
\label{as:R1}
\min_{x\in \R^d} R(x,I_m)=0, \hspace{20 pt} \max_{x\in \R^d} R(x,I_M)=0,
\end{equation}
and there exists constants $K_i>0$ such that, for any $x\in \R^d$, $I\in \R$,
\begin{equation}
\label{as:R2}
-K_1\leq \frac{\partial R}{\partial I}(x,I)\leq-K_1^{-1}<0,
\end{equation}

\begin{equation}
\label{as:R3}
\sup_{\frac{I_m}{2}\leq I \leq 2I_M} \parallel R(\cdot,I)\parallel_{W^{2,\infty}(\R^d)}<K_2.
\end{equation}
Moreover, we make the following assumptions on the initial data:
\beq
\label{as:ue0}
(u_\e^0)_\e \text{ is a sequence of continuous functions which converge in $C_{\rm loc}(\R^d)$ to $u^0$, as $\e \to 0$,}
\eeq
{where $u_\e^0=\e\log n_\e^0$, and there exists a constant $A< \alpha$ and positive constants $C_0$ and $C_1$} such that
\beq
\label{as:u0}
n_\e^0(x)\leq  {\frac{C_0 }{\big(C_1(1+|x|^2)\big)^\f{A}{\e}} },
\eeq
 \beq
 \label{as:I0}
I_m \leq \int_{\R^d} n_\e^0(x)\leq I_M.
\eeq

\subsection {Main results and plan of the paper}

Our main result is the following (see Definition \ref{def:viscosity} for the definition of viscosity sub and supersolutions). 

\begin{theorem}\label{th:main}
Let $n_\e$ be the solution of \fer{SMe} and $u_\e=\e\log n_\e$. Assume \fer{as:R1}--\fer{as:I0}.   Then, along subsequences as $\e\to 0$,  $(I_\e)_\e$ converges a.e. to $I$ and $(u_\e)_\e$ converges locally uniformly to a   function $u$ which is  Lipschitz continuous with respect to $x$ and continuous in $t$, such that
\beq
\label{max}
\|D_x u \|_{L^\infty(\R^d \times \R^+)}\leq 2\alpha,\qquad u(t,x+h)-u(t,x)\leq 2\alpha \log(1+|h|),\qquad \text{for all $x,h\in \R^d$.}
\eeq
Moreover, $u$ is a viscosity supersolution to the following equation
\beq
\label{HJ}
\begin{cases}
 \p_t u -\int_0^\infty \int_{\nu\in S^{d-1}}  \left( e^{k D_x  u\cdot \nu}-1 \right) \f{e^k dSdk}{|e^k-1|^{1+2\al}} =R(x,I),\\
u(x,0)=u^0(x).
\end{cases}
\eeq
For fixed $I$, $u$ is indeed the minimal viscosity supersolution of \fer{HJ} satisfying \fer{max}.
Moreover, $u$ satisfies the following constraint
\beq
\label{constraint}
  \max_{x\in \R^d} u(t,x)=0,\qquad \text{for all $t>0$}.
\eeq
It is also a viscosity subsolution of \fer{HJ} in the following weak sense. {Let} $\vp\in \mathrm{C}^2(\R^+\times \R^d)$ be a test function such that $u-\vp$ takes a maximum at $(t_0,x_0)$ and
\beq
\label{w-test-function}
\vp(t,x+h)-\vp(t,x) \leq (2\alpha-\xi)\log(1+|h|),\qquad \text{for all $(t,x)\in B_r(t_0,x_0)$ and $h\in \R^d$,}
\eeq
with $r$  and $\xi$  positive constants. Then, we have
\beq
\label{w-sub-cr}
\p_t \vp(t_0,x_0) -\int_0^\infty \int_{\nu\in S^{d-1}}  \left( e^{k D_x  \vp(t_0,x_0)\cdot \nu}-1 \right) \f{e^k dSdk}{|e^k-1|^{1+2\al}} \leq \displaystyle\limsup_{  {s\to  t_0}} R(x,I(s)).
\eeq

 \end{theorem}
 
 \noindent
  \noindent 
 A main difficulty in this convergence result is that the Hamiltonian in the above Hamilton-Jacobi equation can take infinite values. Another difficulty comes from the fact that the term $I(t)$ is only BV and potentially discontinuous.
   To prove the convergence of $(u_\e)_\e$ we use the method of semi-relaxed limits \cite{GB.BP:88} in the theory of viscosity solutions. However, since the Hamiltonian in  \fer{HJ} takes infinite values and since the limit $u$ is not in general a viscosity solution of \fer{HJ}, we cannot use this method in a classical manner and further work is required. 
{This issue  is indeed closely related to the work in \cite{CB.EC:13} where a large deviation type result has been obtained for a L\'evy type nonlocal operator where the integral kernel has at most exponential tails. In the case where the integral kernel has exponential tails, a  Hamilton-Jacobi equation close to \fer{HJ}, without the growth term, is obtained at the limit. However, in that case  the function obtained at the limit  does not satisfy necessarily the second regularity result in \fer{max} and it is indeed a viscosity solution to the Hamilton-Jacobi equation.} Note indeed that \fer{max} indicates that there is a strong regularizing effect of the solutions,  independently of the regularity of $R(x,I)$ and the initial condition. Such regularizing effect is proved simultaneously with the proof of the convergence of $(u_\e)_\e$.
 \\

\noindent 
Note that in Theorem \ref{th:main} we do not characterize the limit $u$ as a viscosity solution to a Hamilton-Jacobi equation with constraint, as was the case in the previous results on such selection-mutation models (see for instance \cite{OD.PJ.SM.BP:05,GB.BP:08,GB.SM.BP:09}). We only prove that $u$ is the minimal viscosity supersolution to \fer{HJ} satisfying \fer{max} and a viscosity subsolution in a weak sense. One can wonder if $u$ is indeed a viscosity solution to \fer{HJ}. We do not expect this assertion to be true in general. The fact that the Hamiltonian in \fer{HJ} has infinite values for $|D_x u|\geq 2\alpha$ indicates that \fer{HJ} has a regularizing effect forcing $u$ to verify $|D_x u|\leq 2\alpha$ (i.e. the first property in \fer{max}). However, the second property in \fer{max} is a stronger property and generally is not satisfied by a solution of a Hamilton-Jacobi equation of type \fer{HJ}. In Section \ref{sec:example}, we provide an example of a Hamilton-Jacobi equation of similar type which has a solution that does not satisfy the second inequality in \fer{max}. Existence of such solutions together with the uniqueness of viscosity solutions to Hamilton-Jacobi equations of type \fer{HJ}, with fixed $I$, (see \cite{CB.EC:13}, Section 6) indicates that $u$ might not be a viscosity solution of \fer{HJ} in general. Note that, of course,   thanks to the comparison principle for fixed $I$, $u$ is always greater than (or equal to) the unique viscosity solution of \fer{HJ}. \\

\noindent
The information obtained in Theorem \ref{th:main}   still allows to obtain the concentration of the population's density as Dirac masses, analogously to the previous works \cite{GB.BP:08,GB.SM.BP:09}:
 
 \begin{theorem}
 \label{th:convn}
 Let $n_\e$ be the solution of \fer{SMe}. Assume \fer{as:R1}--\fer{as:I0}. Then, along subsequences as $\e\to 0$, $n_\e$ converges in $\mathrm{L^\infty}\left(\mathrm{w*} (0,\infty) ; \mathcal{M}^1(\R^d) \right)$ to a measure $n$, such that, 
\beq
\label{supp-n}
supp\, n (t,\cdot)\subset \{x \,|\, u(t,x)=0\},\qquad \text{for a.e. $t$}.  
\eeq
Moreover, for all continuous points of $I(t)$, we have
\beq
\label{zero-u-R}
\{x \,|\, u(t,x)=0\}  \subset \{x \,|\, R(x,I(t))=0\} .
\eeq
In particular, if $x\in \R$ and $R$ is monotonic with respect to $x$, then for all $t>0$ except for a countable set of points,
$$
n(t,x) = I(t) \, \delta (x-\overline x(t)).
$$
 \end{theorem}

 \bigskip
 
  \noindent
 The paper is organized as follows. In Section \ref{sec:reg} we provide some preliminary regularity estimates. In section \ref{sec:proof-convergence} we give the main elements of the proof of the convergence of $u_\e$ to a viscosity supersolution of \fer{HJ}. In Section \ref{sec:lem-conv-vB} we prove Proposition \ref{lem:conv-vB} which is an important ingredient in the proof of the convergence of $u_\e$. In sections \ref{sec:proof-thm} and \ref{sec:proof-thm2} we provide respectively the proofs of Theorem \ref{th:main} and Theorem \ref{th:convn}. Finally, in Section \ref{sec:example}, we give an example of a Hamilton-Jacobi equation of type \fer{HJ} which has a viscosity solution not satisfying the second property in \fer{max}.
 \\
 
 \noindent
 Throughout the paper, we denote by $C$ positive constants that are independent of $\e$ but can change from line to line.\\

\section{Regularity estimates}
\label{sec:reg}

In this section we prove the following
\begin{prop}
\label{prop:reg}
Let $(n_\e,I_\e)$ be the solution to \fer{SMe}  and assume \fer{as:R1}--\fer{as:I0}. Then, there exists positive constants $\e_0$ and $C_2$ such that, for all $(t,x)\in \R^+\times \R^d$ and $\e\leq \e_0$,
\beq
\label{boundI}
I_m \leq I_\e(t)\leq I_M,
\eeq
\beq
\label{boundu}
n_\e(t,x)\leq \frac{C_0e^{\f{C_2 t}{\e}}}{ \big(C_1(1+|x|^2)\big)^\f{A}{\e}}.
\eeq
Moreover,  $(I_\e)_\e$ is locally uniformly BV for $\e\leq \e_0$ and hence it converges a.e., as $\e\to 0$ and along subsequences, to a function $I:\R^+\to \R^+$. Moreover, $I$ is nondecreasing in $(0,+\infty)$.
\end{prop}
 
\proof
(i) [Proof of \fer{boundI}] {For $L>L_0$, with $L_0$  a large constant, let $\chi_L$ be} a smooth  function with compact support in $B_L(0)$ such that 
$$
\begin{cases}
\chi_L (x)=1& \text{if $|x|\leq L/2$},\\
0\leq \chi_L (x)\leq 1& \text{if $\f L 2\leq |x|\leq L$},\\
{\|\chi_L\|_{W^{2,\infty}(\R^d)}\leq 1}.
\end{cases}
$${
We define
$$
I_{\e,L}(x)=\int_0^\infty  \int_{\nu\in S^{d-1}}\big(\chi_L(x+(e^{\e k}-1)\nu) - \chi_L(x) \big) \f{e^k }{|e^k-1|^{1+2\al}} dSdk.
$$
It is immediate from the definition of $\chi_L$ that 
\beq
\label{IeL0}
I_{\e,L}(x) \leq 0,\qquad \text{for $|x|\leq \f L2$}.
\eeq
We also prove that there exists a positive constant $C$ such that
\beq
\label{IeLb}
I_{\e,L}(x)\leq C,\qquad \text{for all $\e\leq 1$, $L>L_0$ and $x\in \R^d$}.
\eeq
To this end we split the integral term in $I_{\e,L}$ into two parts:
$$
\begin{array}{rl}
I_{\e,L}&=\int_0^1  \int_{\nu\in S^{d-1}}\big(\chi_L(x+(e^{\e k}-1)\nu) - \chi_L(x) \big) \f{e^k }{|e^k-1|^{1+2\al}} dSdk\\
&+\int_1^\infty  \int_{\nu\in S^{d-1}}\big(\chi_L(x+(e^{\e k}-1)\nu) - \chi_L(x) \big) \f{e^k }{|e^k-1|^{1+2\al}} dSdk\\
&=I_{\e,L,1}+I_{\e,L,2}.
\end{array}
$$
Since $0\leq \chi_L\leq 1$ we have
$$
I_{\e,L,2}\leq \int_1^\infty\int_{\nu\in S^{d-1}} \f{e^k }{|e^k-1|^{1+2\al}} dSdk.
$$
To control $I_{\e,L,1}$ we use the Taylor expansion of $\chi_L(x+(e^{\e k}-1)\nu)$ with respect to $k$ around $k=0$. We compute, for $k\in(0,1)$:
$$
\begin{array}{rl}
\chi_L(x+(e^{\e k}-1)\nu) & = \chi_L(x) +\e k D_x \chi_L(x) \cdot \nu \\
&+\f{k^2}{2}\Big(e^{\e \widetilde k}\e^2 D_x \chi_L \big(y+(e^{\e \widetilde k}-1\big)\nu) \cdot \nu\big)+ e^{2\e \widetilde k}\e^2 \nu^t D^2_{xx} \chi \big(y +(e^{\e \widetilde k}-1)\nu \big) \nu \Big),
\end{array}
$$
with $\widetilde k\in (0,k)$.
We deduce, thanks to the boundedness of the derivatives of $\chi_L$,  that
$$
I_{\e,L,1} \leq \e^2 C\int_{0}^{1}     k^2  \f{e^{(1+2\e)k} }{|e^k-1 |^{1+2\al}} dk,
$$
which is bounded for $\e\leq 1$.
}
\\

\noindent{
We now have at hand a suitable set of test functions that we will use to prove \fer{boundI}.}
We  multiply \fer{SMe} by $\chi_L$ and integrate with respect to $x$ and obtain, using  Fubini's Theorem, {\fer{IeL0} and \fer{IeLb}},
{$$
\begin{array}{rl}
\e\f{d}{dt} \int_{\R^d} \chi_L(x) n_\e(t,x)dx & = \int_{\R^d}  n_\e(t,x)\left( \int_0^\infty  \int_{\nu\in S^{d-1}}\big(\chi_L(x+(e^{\e k}-1)\nu) - \chi_L(x) \big) \f{e^k }{|e^k-1|^{1+2\al}} dSdk\right) dx
\\
&+ \int_{\R^d} \chi_L(x) n_\e(t,x)R(x,I_\e(t))dx\\
& \leq C\int_{\f L2}^\infty n_\e(t,x)dx+ \int_{\R^d} \chi_L(x) n_\e(t,x)R(x,I_\e(t))dx.
\end{array}
$$}
We then let $L$ go to $+\infty$, {and use the fact that $n_\e(t,\cdot )\in L^1(\R^d)$} to obtain
$$
\e\f{d}{dt}I_\e(t) =\int_{\R^d} n_\e(t,x)R(x,I_\e(t))dx.
$$
Using the above equation, \fer{as:R1} and \fer{as:I0} we obtain \fer{boundI}.\\

\noindent
(ii) [Proof of \fer{boundu}] We define, for $C_2$ a positive constant,
$$
s(t,x)= \frac{C_0e^{ \f{C_2 t}{\e}}}{\big(C_1(1+|x|^2)\big)^\f{A}{\e}}.
$$
We show that, for $C_2$ large enough, $s$ is a supersolution to \fer{SMe}. Note that \fer{SMe}, with $I_\e$ fixed, admits a comparison principle, since \fer{SM} admits a comparison principle (see \cite{GB.CI:08}--Theorem 3). Moreover, thanks to Assumption \fer{as:u0},
$$
n_\e^0(x)\leq  \frac{C_0 }{\big(C_1(1+|x|^2)\big)^\f{A}{\e}} = s(0,x).
$$
We hence obtain \fer{boundu} thanks to the comparison principle.\\

\noindent
To prove that, for $C_2$ large enough, $s$ is a supersolution to \fer{SMe}, since $R(x,I_\e)$ is bounded thanks to \fer{as:R3}, it is enough to prove that
$$
F(t,x)=\int_0^\infty \int_{\nu\in S^{d-1}} \left( s(t,x+(e^{\e k}-1)\nu) - s(t,x) \right) \f{e^k }{|e^k-1|^{1+2\al}}dSdk\leq C s(t,x),
$$
for $C$ a constant which is large enough but is independent of $\e$.
We compute
$$
\f{F(t,x)}{s(t,x)}=\int_0^\infty \int_{\nu\in S^{d-1}}  \left( \f{(1+|x|^2)^{\f{A}{\e}}}{(1+|x+(e^{\e k}-1)\nu|^2)^{\f{A}{\e}}} - 1 \right) \f{e^k }{|e^k-1|^{1+2\al}}dSdk.
$$
We split the above integral into two parts, that we will control separately,
$$
G_1=\int_1^\infty \int_{\nu\in S^{d-1}}  \left( \f{(1+|x|^2)^{\f{A}{\e}}}{(1+|x+(e^{\e k}-1)\nu|^2)^{\f{A}{\e}}} - 1 \right) \f{e^k }{|e^k-1|^{1+2\al}}dSdk,
$$
$$
G_2=\int_0^1 \int_{\nu\in S^{d-1}}  \left( \f{(1+|x|^2)^{\f{A}{\e}}}{(1+|x+(e^{\e k}-1)\nu|^2)^{\f{A}{\e}}} - 1 \right) \f{e^k }{|e^k-1|^{1+2\al}}dSdk.
$$
In order to control the above integrals we use the following inequality:
$$
\f{1+|y-l|^2}{1+|y|^2}\leq (1+|l|)^2, \quad \text{with} \quad l=(e^{\e k}-1)\nu, \quad x=y-l.
$$
We deduce that
\beq
\label{bounds}
 s(x,k,\nu) := \f{(1+|x|^2)^{\f{A}{\e}}}{(1+|x+(e^{\e k}-1)\nu|^2)^{\f{A}{\e}}} \leq e^{2Ak},
\eeq
 and hence
 $$
G_1\leq \int_1^\infty \int_{\nu\in S^{d-1}}  (e^{2Ak}-1)\f{e^{k} }{|e^k-1|^{1+2\al}}dSdk \leq C. 
 $$
 Note that the above integral is bounded since $A<\alpha$.\\
 
 \noindent
 In order to control $G_2$ we use the Taylor's expansion of  $ s(x,k,\nu)$ with respect to $k$, around $k=0$. We compute
 $$
 \f{\p}{\p k} s(x,k,\nu)= -2A \f{(1+|x|^2)^{\f{A}{\e}}}{(1+|x+(e^{\e k}-1)\nu|^2)^{1+\f{A}{\e}}} 
 {(x+(e^{\e k}-1)\nu,\nu)}e^{\e k},
 $$
  $$
  \begin{array}{rl}
 \f{\p^2}{\p k^2} s(x,k,\nu) &= 4A(A+\e) \f{(1+|x|^2)^{\f{A}{\e}}}{(1+|x+(e^{\e k}-1)\nu|^2)^{2+\f{A}{\e}}} {(x+(e^{\e k}-1)\nu,\nu)}^2e^{2\e k}\\
 &-2A\e\, \f{(1+|x|^2)^{\f{A}{\e}}}{(1+|x+(e^{\e k}-1)\nu|^2)^{1+\f{A}{\e}}} {(x+(e^{\e k}-1)\nu,\nu)}e^{\e k}\\
 &-2A\e\, \f{(1+|x|^2)^{\f{A}{\e}}}{(1+|x+(e^{\e k}-1)\nu|^2)^{1+\f{A}{\e}}} {(\nu,\nu)}e^{2\e k}.
 \end{array}
 $$
One can verify that
$$
 \int_0^1 \int_{\nu\in S^{d-1}}  k\f{\p}{\p k} s(x,0,\nu) \f{e^k }{|e^k-1|^{1+2\al}}dSdk=0,
$$
and  thanks to \fer{bounds},  $\f{\p^2}{\p k^2} s(x,k,\nu)$ is bounded for  $\e\leq \e_0$ and $0\leq k\leq 1$. We deduce that, for all $\e\leq \e_0$,
 $$
G_2\leq \widetilde C\int_0^1 \int_{\nu\in S^{d-1}}  k^2\f{e^{k} }{|e^k-1|^{1+2\al}}dSdk \leq C.
 $$
Combining the above inequalities on $G_1$ and $G_2$ we obtain that $F(t,x)\leq Cs(t,x)$ for $C$ large enough and $\e\leq \e_0$.\\

\noindent 
(iii) [Uniform BV bound on $I_\e$]
The proof of uniform BV bound on $I_\e$ is an adaptation of the proof of Theorem 3.1. in \cite{GB.BP:08}.
Integrating \fer{SMe} with respect to $x$ we obtain
$$
\f{d}{dt}I_\e=\f{1}{\e} \int_{\R^d} n_\e(t,x) \,R(x,I_\e(t)) dx.
$$
We define $ J_\e(t)=\f{1}{\e} \int_{\R^d} n_\e(t,x) \,R(x,I_\e(t)) dx$ and  $ J_{\e,L}(t)=\f{1}{\e} \int_{\R^d} \chi_L(x) n_\e(t,x) \,R(x,I_\e(t)) dx$, with $\chi_L$ defined in the proof of part (i). We then differentiate $J_{\e,L}$ with respect to $t$ and obtain
$$
\begin{array}{rl}
\f{d}{dt} J_{\e,L}(t) &= \f{1}{\e} \displaystyle\int_{\R^d}  \chi_L(x)\f{\p}{\p t} n_\e(t,x) \,R(x,I_\e(t)) dx +\f{1}{\e}\Big(\displaystyle \int_{\R^d} \chi_L(x) n_\e(t,x) \,\f{\p}{\p I}R(x,I_\e(t)) dx \Big) \f{d}{dt} I_\e(t)\\
 &=   \displaystyle \f{1}{\e^2}\int_{\R^d}\int_0^\infty \int_{\nu\in S^{d-1}} \chi_L(x)\left( n_\e(t,x+(e^{\e k}-1)\nu) - n_\e(t,x) \right) R(x,I_\e(t))\f{e^k }{|e^k-1|^{1+2\al}}dSdk dx\\
 &+\displaystyle \f{1}{\e^2} \int_{\R^d} {\chi_L(x)} n_\e(t,x) \,R(x,I_\e(t))^2 dx +\f{1}{\e}\Big(\displaystyle \int_{\R^d}\chi_L(x) n_\e(t,x) \,\f{\p}{\p I}R(x,I_\e(t)) dx \Big) \f{d}{dt} I_\e(t)\\
 &=A_1+A_2+A_3.
\end{array}
$$ 
We rewrite $A_1$ as below
$$
\begin{array}{l}
A_1=\\   \f{1}{\e^2}\int_{\R^d}\int_0^\infty \int_{  S^{d-1}} n_\e \big(t,x+(e^{\e k}-1)\nu\big)  \left( (\chi_L R)(x,I_\e ) -(\chi_L R) \big(x+(e^{\e k}-1)\nu,I_\e )\big)\right) \f{e^k }{|e^k-1|^{1+2\al}}dSdk dx  +\\ 
  \f{1}{\e^2}\int_{\R^d}\int_0^\infty \int_{  S^{d-1}}   \left( (n_\e\chi_L R) \big(x+(e^{\e k}-1)\nu,I_\e )\big)-    (n_\e \chi_L R) (x,I_\e ) \right)\f{e^k }{|e^k-1|^{1+2\al}} dSdkdx =
\\
 \f{1}{\e^2}\int_0^\infty \int_{  S^{d-1}}  \int_{\R^d} n_\e \big(t,x+(e^{\e k}-1)\nu\big)  \left( (\chi_L R)(x,I_\e ) -(\chi_L R) \big(x+(e^{\e k}-1)\nu,I_\e )\big)\right) \f{e^k }{|e^k-1|^{1+2\al}}dx dSdk  =\\
  \f{1}{\e^2} \int_0^\infty \int_{  S^{d-1}}\int_{\R^d} n_\e \big(t,y\big)  \left( (\chi_L R)(y-(e^{\e k}-1)\nu,I_\e ) -(\chi_L R) \big(y,I_\e )\big)\right) \f{e^k }{|e^k-1|^{1+2\al}}dy dSdk .
\end{array}
$$ 
Note that here we have used  Fubini's theorem  on the first and the second integral term in order to integrate with respect to $x$ before integrating with respect to $\nu$ and then $k$, which allows to show in particular that the second integral  term is null. We then use \fer{as:R3} and a Taylor expansion of the integrand of the last line with respect to $\e$ around $\e=0$ to obtain that, for $\e\leq \e_0$ small enough, there exists a positive constant $C$, independent of $\e$ and $L$, such that
$$
|A_1| \leq C.
$$
We next notice that $A_2$ is positive. We hence obtain that 
$$
\f{d}{dt} J_{\e,L}(t) \geq -C+ \f{1}{\e}\Big(\displaystyle \int_{\R^d}\chi_L(x) n_\e(t,x) \,\f{\p}{\p I}R(x,I_\e(t)) dx \Big) J_\e(t).
$$
We then let $L$ go to $+\infty$ and use \fer{boundu} and \fer{as:R2} to obtain
$$
\f{d}{dt} J_{\e}(t) \geq -C+ \f{1}{\e}\Big(\displaystyle \int_{\R^d}  n_\e(t,x) \,\f{\p}{\p I}R(x,I_\e(t)) dx \Big) J_\e(t).
$$
It follows, thanks to \fer{as:R2} and \fer{boundI}, that for $\e\leq \e_0$,
$$
\f{d}{dt} (J_\e(t))_- \leq C- \f{I_m }{\e K_1} (J_\e(t))_-,
$$
with $(J_\e(t))_-=\max(0,-J_\e(t))$. We deduce that
$$
(J_\e(t))_- \leq \e\f{C  K_1}{I_m} + (J_\e(0))_- e^{-\f{I_m  t}{ K_1\e}}.
$$
We next use \fer{as:R3} and \fer{as:I0} to obtain
$$
(J_\e(0))_-\leq \f{1}{\e} \left(\int_{\R^d} n_\e^0(x) \,R(x,I_\e(0)) dx\right)_- \leq \f{K_2I_M}{\e}.
$$
  We deduce that
\beq
\label{Je-}
(J_\e(t))_- \leq \e\f{C K_1}{I_m} + \f{K_2I_M}{\e} e^{-\f{I_m t}{K_1\e}}.
\eeq
Finally, we show that the above inequality leads to a BV estimate on $I_\e$. To this end, we compute using \fer{boundI} and the above inequality:
$$
\begin{array}{rl}
\displaystyle \int_0^T | \f{d}{dt} I_\e(t) |dt &=\displaystyle  \int_0^T  \f{d}{dt} I_\e(t) dt +2\int_0^T ( J_\e(t))_- dt \\
&\displaystyle  \leq I_M-I_m+2 \e\f{CK_1T}{I_m }+\f{K_2K_1I_M}{ I_m}.
\end{array}
$$
We conclude that $(I_\e)_\e$ is locally uniformly BV for $\e\leq \e_0$. As a consequence, $(I_\e)_\e$ converges a.e., as $\e\to 0$ and along subsequences, to a function $I:\R^+\to \R^+$. Moreover, for all $t_0>0$, $I$ is nondecreasing in $[t_0,+\infty)$ thanks to \fer{Je-}.
\qed

\section{Convergence of $u_\e$ to a viscosity supersolution of \fer{HJ}}
\label{sec:proof-convergence}

In this Section, we prove the following 
\begin{prop}
\label{prop:conv}
Assume \fer{as:R1}--\fer{as:I0}. As $\e\to 0$ and along subsequences, $u_\e$ converges to $u$, a viscosity supersolution of the Hamilton-Jacobi equation in \fer{HJ}. Moreover, $u$ satisfies \fer{max}.
\end{prop}

\proof
For a technical reason, we will need to deal with an equation with negative growth rate. Therefore, we modify $n_\e$ in the following way
$$
m_\e(t,x)=n_\e(t,x)\, e^{-\f{K_2t}{\e}}.
$$
The above function solves
\beq
\label{eq:me}
\e \p_t m_\e (t,x) = \int_0^\infty \int_{\nu\in S^{d-1}} \left( m_\e(t,x+(e^{\e k}-1)\nu) - m_\e(t,x) \right) \f{e^k }{|e^k-1|^{1+2\al}}dSdk+m_\e(t,x) \,{\widetilde R(x,I_\e(t))},
\eeq
with 
\beq
\label{Gneg}
{\widetilde R(x,I)=R(x,I)-K_2} \leq 0. 
\eeq
We then define 
$$
v_\e = \e \log (m_\e).
$$
It is easy to verify that $(u_\e)_\e$ converges to $u$ a viscosity supersolution to the Hamilton-Jacobi equation in \fer{HJ} if and only if $(v_\e)_\e$ converges to $v$, a viscosity supersolution of the following equation
\beq
\label{HJ-G}
\p_t v -H( D_x  v)  ={\widetilde R(x,I)},
\eeq
with 
\beq
\label{def:H}
H(D_x v)=\int_0^\infty \int_{\nu\in S^{d-1}}  \left( e^{k D_x  v\cdot \nu}-1 \right) \f{e^k dSdk}{|e^k-1|^{1+2\al}}.
\eeq
\noindent
{Moreover, if $v$ satisfies \fer{max}, then $u$ also satisfies \fer{max}.}\\

\noindent
In what follows we will prove that $(v_\e)_\e$ converges indeed to a viscosity supersolution of \fer{HJ-G} {that satisfies \fer{max}}. 
We first notice that thanks to \fer{boundu}, $v_\e$ is locally uniformly bounded {above by}:
$$
v_\e(t,x) \leq \e\log C_0 -A \log(C_1)+(C_2-K_2) t -A \log (1+|x|^2).
$$
To avoid lower estimate we use a classical trick by modifying $v_\e$ a little bit ({see for instance \cite{GB.BP:90}}):
\beq
\label{def:veB}
v_\e^B=\e \log (m_\e+e^{-\f{B}{\e}}),
\eeq
with $B$ a large positive constant.
One can verify that $v_\e^B$ is locally uniformly bounded from above and below.\\
We prove the following results:
\begin{prop}
\label{lem:conv-vB}
Assume that $(I_\e)_\e$ converges as $\e\to 0$ to $I$. Then, as $\e \to 0$, the sequence $(v_\e^B)_\e$ converges to $v^B$ a viscosity supersolution to the following equation  
\beq
\label{eq:vB}
\min(v^B+B, \p_t v^B -H( D_x  v^B)  -{ \widetilde R(x,I)})=0.
\eeq
Moreover, $v^B$ satisfies
\beq
\label{LipvB}
\|D_x v^B \|_{L^\infty(\R^d \times \R^+)}\leq 2\alpha,\qquad   v^B(t,x+h)- v^B(t,x)\leq 2\alpha \log(1+|h|),
\eeq
for all $(x, h) \in \R^d\times \R^d$.
\end{prop}
\begin{lemma}
\label{lem:vBv}
For any compact set $K\in \R^+\times \R^d$, there exists $B_0$ large enough, such that for all $B\geq B_0$, $-B <v^B$ in $K$.
\end{lemma}
We postpone the proof of Proposition \ref{lem:conv-vB} to the next section and the proof of Lemma \ref{lem:vBv} to the end of this paragraph and explain first how they allow us to conclude.\\
Thanks to Proposition \ref{prop:reg}, as $\e\to 0$ and along   subsequences, $I_\e$ converges a.e. to a function $I$. Proposition \ref{lem:conv-vB} implies that, along such subsequences and for all $B>0$, $v_\e^B$ converges to $v^B$ a viscosity supersolution of \fer{eq:vB}.
Let's  fix a compact set $K$ and consider $B_0$ given by Lemma \ref{lem:vBv}. Thanks to the definitions of $v_\e$ and $v_\e^B$ we can write
$$
v_\e=\e \log (e^{\f{v_\e^{B_0}}{\e}}-e^{\f{-B_0}{\e}} ).
$$
 We then   use the fact that, in the set $K$, $-B_0 <v^{B_0}$, to obtain that $v_\e$ converges, in the set $K$, to $v=v^{B_0}$. Moreover, $v$ is a viscosity supersolution to  \fer{HJ-G}, in the set $K$  thanks to \fer{eq:vB} and it satisfies \fer{max} thanks to \fer{LipvB}.
 \qed

 \bigskip
\noindent
To prove Lemma  \ref{lem:vBv} we first introduce the following semi-relaxed limits
$$
\overline{u}(t,x)=\limsup_{\underset{\e\to 0}{(s,y)\to (t,x)}} u_\e (s,y),\qquad {\overline{v}(t,x)=\limsup_{\underset{\e\to 0}{(s,y)\to (t,x)}} v_\e (s,y),}
$$
$$
\overline{v}^B(t,x)=\limsup_{\underset{\e\to 0}{(s,y)\to (t,x)}} v_\e^B (t,x),
 \qquad
\underline{v}^B(t,x)=\liminf_{\underset{\e\to 0}{(s,y)\to (t,x)}} v_\e^B (t,x).
$$
Note that we can define such quantities, since $u_\e$ is locally uniformly bounded from above and $v_\e^B$ is locally uniformly bounded from below and above.
We then  prove the following lemma.
\begin{lemma}
\label{lem:maxubar}
Assume \fer{as:R1}--\fer{as:I0}. Then, for all $t \in \R^+$, we have
$$
\max_{x\in \R^d} \; \overline u(t,x)\geq 0.
$$
\end{lemma}

{\bf [Proof of Lemma \ref{lem:maxubar}]} Let's fix $t\in \R^+$ and assume that $\max_{x\in \R^d} \; \overline u(t,x)= -a<0$. Note that such maximum is attained thanks to \fer{boundu}. Thanks to \fer{boundu} there exists  constants $r>0$  large enough and $\e_0$ small enough such that, for all $\e\leq \e_0$,
$$
\int_{B_r(0)^c} n_\e(t,x) dx <\f{I_m}{2}.
$$
It follows from \fer{boundI} that
$$
\f{I_m}{2 } \leq \int_{B_r(0)} n_\e(t,x) dx =  \int_{B_r(0)} e^{\f{u_\e(t,x)}{\e}} dx.
$$ 
Letting $\e\to 0$ we obtain that
$$
\f{I_m}{2 } \leq    \int_{B_r(0)} \limsup_{\e\to 0} e^{\f{u_\e(t,x)}{\e}} dx =0,
$$ 
 which is a contradiction.
\qed

\bigskip

{\bf [Proof of Lemma \ref{lem:vBv}]}
Let's fix $T>0$. Thanks to Lemma \ref{lem:maxubar} and the definition of $\overline v^B$ we have
$$
-K_2T \leq \max_{x} v^B(t,x) = \max_{x}  \overline  v^B(t,\cdot) ,\qquad \text{for $t\in [0,T] $}.
$$
We also note that, thanks to \fer{boundu}, there exists a positive constant $r$ large enough such that
$$
\overline v(t,x) = \overline u(t,x)-K_2{t} < -K_2T, \qquad \text{for  $(t,x)\in [0,T]\times B_r(0)^c$}.
$$
It follows that for $B>K_2 T$ and $t\in[0,T]$, $\overline v^B(t,\cdot)$ attains its maximum with respect to $x$ in the set $B_r(0)$. Moreover, this maximum is greater than $-K_2T$. 
Next, using the Lipschitz continuity of $v^B$ given by Proposition \ref{lem:conv-vB}, we deduce that for any compact set $K\subset [0,T]\times \R^d$ and $B>K_2T$, there exists a constant $C$ large enough, independent of $B$, such that
$$
-C< v^B(t,x) ,\qquad \text{for all $(t,x)\in K$}.
$$
Finally, taking $B_0=\max( K_2T+1,C )$ we conclude that
$$
-B_0 < v^{B_0}(t,x),\qquad \text{for all $(t,x)\in K$}.
$$
\qed

 \section{Proof of Proposition \ref{lem:conv-vB}}
 \label{sec:lem-conv-vB}
 
To prove Proposition \ref{lem:conv-vB}, we will work with semi-relaxed limits $\overline{v}^B$ and $\underline{v}^B$.
A classical method in the theory of viscosity solutions is to prove that 
$\overline{v}^B$ and $\underline{v}^B$ are respectively sub and supersolutions of \fer{eq:vB} and then use a comparison principle to obtain that $\overline{v}^B\leq  \underline{v}^B$. This would imply that $\overline{v}^B=  \underline{v}^B$ and that $(v_\e)_\e$ converges locally uniformly to the solution of \fer{eq:vB}. Here, we cannot use this strategy because $\overline{v}^B$ is not generally a subsolution of \fer{eq:vB}. To overcome this difficulty we first regularize the supersolution $\underline{v}^B$ and modify it to become a strict supersolution and to satisfy some required properties. Then we use it as a test function that we compare with $\overline{v}^B$ to obtain directly that $\overline{v}^B \leq \underline{v}^B$. See \cite{CB.EC:13,GB.AB.EC:14} where this method has been suggested in other contexts.\\

\noindent
Before providing the proof of  Proposition  \ref{lem:conv-vB} we first recall the definition of viscosity solutions for \fer{HJ} which has a discontinuous Hamiltonian (see \cite{GB:94}-page 80). {Note that here the discontinuity of the Hamilton-Jacobi equation comes from the fact that the function $I$ can be potentially discontinuous since it is only of bounded variation.}

\begin{definition}[viscosity solutions]
\label{def:viscosity}
 
(i) An upper semi-continuous function $u$ which is locally bounded is a subsolution of \fer{HJ} if and only if 
$$
\forall \vp \in \mathrm{C}^2(\R^+\times \R^d), \text{if $u-\vp$ takes a local maximum at $(t_0,x_0)$, then}
$$
$$
\p_t \vp(t_0,x_0) -\int_0^\infty \int_{\nu\in S^{d-1}}  \left( e^{k D_x  \vp(t_0,x_0)\cdot \nu}-1 \right) \f{e^k dSdk}{|e^k-1|^{1+2\al}} \leq \displaystyle\limsup_{  {s\to  t_0}} R(x,I(s)).
$$
(ii) A lower semi-continuous function $u$ which is locally bounded is a supersolution of \fer{HJ} if and only if 
$$
\forall \vp \in \mathrm{C}^2(\R^+\times \R^d), \text{if $u-\vp$ takes a local minimum at $(t_0,x_0)$, then}
$$
$$
\p_t \vp(t_0,x_0) -\int_0^\infty \int_{\nu\in S^{d-1}}  \left( e^{k D_x  \vp(t_0,x_0)\cdot \nu}-1 \right) \f{e^k dSdk}{|e^k-1|^{1+2\al}} \geq \displaystyle\liminf_{  {s\to  t_0}} R(x,I(s)).
$$
(iii) A continuous function $u$ which satisfies {both the} properties above, is a viscosity solution of \fer{HJ}.
\end{definition}

\bigskip

\noindent
We provide the proof of Proposition \ref{lem:conv-vB} in several steps. 
Note first that replacing   \fer{def:veB} in \fer{eq:me}, we obtain that
\beq
\label{eq:veB}
\p_t v_\e^B=\int_0^\infty \int_{\nu \in S^{d-1}} (e^{\f{v_\e^B(t,x+(e^{\e k}-1)\nu) - v_\e^B(t,x)}{\e}} -1) \f{e^k}{|e^k-1|^{1+2\alpha}}dS dk+\f{m_\e}{m_\e+e^{\f{-B}{\e}}} {\widetilde R(x,I_\e(t))}.
\eeq

(i) We first prove that $\underline{v}^B(t,x)$ is a supersolution to \fer{eq:vB}. \\
To this end, let's suppose that $\vp\in \mathrm{C}( \R^+\times \R^d;\R) \cap \mathrm{C}^2(\mathcal O(t_0,x_0))$, with $\mathcal O(t_0,x_0)$ a neighborhood of $(t_0,x_0)$, is a test function such that $\underline{v}^B-\vp$ attains a  global and strict minimum at $(t_0,x_0)$. Then, (see \cite{GB:94}, Lemma 4.2) there exists a sequence $(t_\e,x_\e)_\e$ such that $v_\e^B -\vp $ has a global minimum at $(t_\e,x_\e)$, $(t_\e,x_\e) \to (t_0,x_0)$ as $\e\to 0$ and $v_\e^B(t_\e,x_\e)\to \underline v^B (t_0,x_0)$. \\
Using \fer{eq:veB}, \fer{Gneg} and the fact that $(t_\e,x_\e)$ is a global minimum point of $v_\e^B-\vp$, we obtain
\beq
\label{eq:supertest}
\begin{array}{rl}
\p_t \vp(t_\e,x_\e) &\geq \int_0^M \int_{\nu \in S^{d-1}} (e^{\f{\vp(t_\e,x_\e+(e^{\e k}-1)\nu) - \vp(t_\e,x_\e)}{\e}} -1)\f{e^k}{|e^k-1|^{1+2\alpha}}dS dk+{\widetilde R(x_\e,I_\e(t_\e))}\\
&+\int_M^\infty \int_{\nu \in S^{d-1}} (e^{\f{v_\e^B(t_\e,x_\e+(e^{\e k}-1)\nu) - v_\e^B(t_\e,x_\e)}{\e}} -1)\f{e^k}{|e^k-1|^{1+2\alpha}}dS dk.
\end{array}
\eeq
{Note that the above formula holds for all $M>0$ since
$$
\vp(t_\e,z) - \vp(t_\e,x_\e) \leq v_\e^B(t_\e,z) - v_\e^B(t_\e,x_\e),\qquad \text{for all $z\in \R^d$.}
$$}
Using the Taylor-Lagrange formula we have, for $\e$ small enough and thanks to the fact that $\vp\in C^2(\mathcal O(t_0,x_0))$, for $\mu \in (0,\e)$,
$$
\begin{array}{rl}
\vp(t_\e,x_\e+(e^{\e k}-1)\nu) & = \vp(t_\e,x_\e) +k\e D_x \vp(t_\e,x_\e) \cdot \nu \\
&+\f{\e^2}{2}\Big(e^{\mu k}k^2 D_x \vp \big(t_\e,x_\e +(e^{\mu k}-1\big)\nu) \cdot \nu\big)+ e^{2\mu k}k^2 \nu^t D^2_{xx} \vp \big(t_\e,x_\e +(e^{\mu k}-1)\nu \big) \nu \Big).
\end{array}
$$
Therefore, for fixed $M$ as $\e \to 0$, the first integral term at the {r.h.s.} of \fer{eq:supertest} converges to 
$$
 \int_0^M \int_{\nu \in S^{d-1}} (e^{k D_x  \vp(t_0,x_0) \cdot \nu } -1)\f{e^k}{|e^k-1|^{1+2\alpha}}dS dk.
 $$
Passing to the limit in  \fer{eq:supertest} as $\e\to 0$ we thus obtain that
$$
\begin{array}{rl}
\p_t \vp(t_0,x_0) &\geq \int_0^M \int_{\nu \in S^{d-1}}(e^{k D_x  \vp(t_0,x_0) \cdot \nu } -1) \f{e^k}{|e^k-1|^{1+2\alpha}}dS dk+\displaystyle\liminf_{\underset{\e\to 0}{(s,y)\to (t_0,x_0)}} {\widetilde R(y,  I_\e(s))},\\
&-\int_M^\infty \int_{\nu \in S^{d-1}}  \f{e^k}{|e^k-1|^{1+2\alpha}}dS dk,
\end{array}
$$
where we have used the positivity of the exponential term in the last term of \fer{eq:supertest}. Letting $M\to \infty$ and using the smoothness of $\widetilde R$ with respect to the first variable and it's monotonicity with respect to its second variable we obtain that
$$ 
\p_t \vp(t_0,x_0)  \geq \int_0^{+\infty} \int_{\nu \in S^{d-1}}(e^{k D_x  \vp(t_0,x_0) \cdot \nu } -1) \f{e^k}{|e^k-1|^{1+2\alpha}}dS dk+   {\widetilde R(x_0,\limsup_{\underset{\e\to 0}{s\to t_0}}  I_\e(s))},
$$
To prove that $\underline{v}^B(t,x)$ is a supersolution to \fer{eq:vB}, that is
$$ 
\p_t \vp(t_0,x_0)  \geq \int_0^{+\infty} \int_{\nu \in S^{d-1}}(e^{k D_x  \vp(t_0,x_0) \cdot \nu } -1) \f{e^k}{|e^k-1|^{1+2\alpha}}dS dk+ {  \widetilde R(x_0,\overline I(t_0) )},
$$
with
$$
{ \overline I(t_0)=\limsup_{s\to t_0} I(s)},
$$
it remains to prove that
$$
\limsup_{\underset{\e\to 0}{s\to t_0}}  {I_\e(s)} \leq \limsup_{s\to t_0} I(s).
$$
This can be proved similarly to the proof of Theorem 4.1. in \cite{GB.BP:08}.\\
Finally, from \fer{def:veB} it is immediate that $\underline v^B\geq -B$. Therefore, $\underline{v}^B(t,x)$ is a supersolution to \fer{eq:vB}.
\\

(ii) We prove that
\beq
\label{uinf-lip}
\|D_x \underline v^B \|_{L^\infty(\R^d \times \R^+)}\leq 2\alpha,\qquad \underline v^B(t,x+h)-\underline v^B(t,x)\leq 2\alpha \log(1+|h|).
\eeq
{Let $(\overline t,\overline x)\in \R^+\times \R^d$ be such that there exists a test function $\vp\in \mathrm{C}( \R^+\times \R^d ; \R) \cap \mathrm{C}^2(\mathcal O(\overline t,\overline x))$, with $\mathcal O(\overline t,\overline x)$ a neighborhood of $(\overline t,\overline x)$,  such that $\underline{v}^B-\vp$ attains a  global and strict minimum at $(\overline t,\overline x)$. We first prove the second inequality of \fer{uinf-lip} for such points. Note that the set of such points is dense in $\R^+\times \R^d$.}\\

\noindent
Let's suppose that there exist $(k_0,\nu_0)\in \R^+\times S^{d-1}$ and $b>0$, such that 
$$
\underline v^B (\overline t , \overline x+(e^{k_0}-1)\nu_0) -\underline v^B(\overline t,\overline x) \geq 2k_0\alpha +b.
$$
Since $\underline v^B$ is lower semi-continuous, we deduce that there exist positive constants $k_1$ and $k_2$ such that $k_1 <k_0<k_2$ and an open set $\Omega_0\subset S^{d-1}$ such that $\nu_0\in \Omega_0$ and
$$
\underline v^B (\overline t,\overline x+(e^{k}-1)\nu) -\underline v^B(\overline t,\overline x) \geq 2k\alpha +\f{b}{2}, \quad \text{for $k\in (k_1,k_2)$, and $\nu \in \Omega_0$}.
$$
From the definition of $\underline v^B$, we also deduce that, there exists a subsequence $(\e_n)_n$, with $\e_n\to 0$ as $n\to +\infty$,   and there exists $(t_n,x_n)$ such that {$(t_n,x_n)$ is a global minimum point of $v_{\e_n}^B-\vp$,} $(t_n,x_n)\to (\overline t,\overline x)$  and $v_{\e_n}^B(t_n,x_n)\to  \underline v^B(\overline t,\overline x)$, as $n\to +\infty$, and 
$$
 v_{\e_n}^B (t_n,x_n+(e^{k}-1)\nu) -  v_{\e_n}^B(t_n,x_n) \geq 2k\alpha +\f{b}{4}, \quad \text{for $k\in (k_1,k_2)$, and $\nu \in \Omega_0$}.
$$
{Similarly to} \fer{eq:supertest} we obtain that
$$
\begin{array}{rl}
\p_t \vp(t_n,x_n) &\geq \int_0^M \int_{\nu \in S^{d-1}} (e^{\f{\vp(t_n,x_n+(e^{\e_n k}-1)\nu) - \vp(t_n,x_n)}{\e_n}} -1)\f{e^k}{|e^k-1|^{+2\alpha}}dS dk+{ \widetilde R(x_n,I_{\e_n}(t_n))}\\
&+\int_{\f{k_1}{\e_n}}^{\f{k_2}{\e_n}} \int_{\nu \in \Omega_0}  e^{\f{v_{\e_n}^B(t_n,x_n+(e^{\e_n k}-1)\nu) - v_{\e_n}^B(t_n,x_n)}{\e_n}} \f{e^k}{|e^k-1|^{1+2\alpha}}dS dk
-\int_M^\infty \int_{\nu \in S^{d-1}}  \f{e^k}{|e^k-1|^{1+2\alpha}}dS dk\\
 &\geq \int_0^M \int_{\nu \in S^{d-1}} (e^{\f{\vp(t_n,x_n+(e^{\e_n k}-1)\nu) - \vp(t_n,x_n)}{\e_n}} -1)\f{e^k}{|e^k-1|^{1+2\alpha}}dS dk+{\widetilde R(x_n,I_{\e_n}(t_n))}\\
&+\int_{k_1 }^{k_2} \int_{\nu \in \Omega_0} e^{\f{2k\alpha+\f{b}{4}}{\e_n}} \f{e^{\f{k}{\e_n}}}{e^{\f{k(1+2\alpha)}{\e_n}}}dS  \f{dk}{\e_n}
-\int_M^\infty \int_{\nu \in S^{d-1}}  \f{e^k}{|e^k-1|^{1+2\alpha}}dS dk
\end{array}
$$
Note that the third term in the r.h.s. of the above inequality goes to $+\infty$ as $n\to +\infty$, while the other terms are bounded and asymptotically, as $n\to +\infty$, greater than 
$$
 \int_0^M \int_{\nu \in S^{d-1}}(e^{k D_x  \vp(\overline t,\overline x) \cdot \nu } -1) \f{e^k}{|e^k-1|^{1+2\alpha}}dS dk+{\widetilde R(\overline x, \overline I(\overline t))}
 -\int_M^\infty \int_{\nu \in S^{d-1}}  \f{e^k}{|e^k-1|^{1+2\alpha}}dS dk.
$$
This is in contradiction with the fact that $\p_t \vp(t_n,x_n)$ is bounded, and hence we obtain the second inequality in \fer{uinf-lip}:
\beq
\label{log-xbar}
 \underline v^B({\overline t} ,{\overline x}+h)-\underline v^B({\overline t},{\overline x})\leq 2\alpha \log(1+|h|), \quad \text{for all $h\in \R^d$}.
\eeq
  {We hence have proved \fer{uinf-lip} for all $(\overline t,\overline x)$ below which we can put a $C^2$ test function. Note also that since 
  $(\overline t, \overline x)$ is a global minimum point of $\underline v^B-\vp$, we have
$$
\vp({\overline t} ,{\overline x}+h)-\vp({\overline t},{\overline x})\leq 2\alpha \log(1+|h|), \quad \text{for all $h\in \R^d$}.
$$
and hence
$$
|\nabla \vp|(\overline t,\overline x) \leq 2\al.
$$
We deduce that
$$
-|\nabla \underline  v^B| \geq -2\al,\qquad \text{ in $\R^+\times \R^d$},
$$
in the viscosity sense. As a consequence, $\underline v^B$ is Lipschitz continuous with Lipschitz constant $2\al$. Since the set of the points $(\overline t,\overline x)$ below which we can put a $C^2$ test function is dense in $ \R^+\times \R^d$, the continuity of $\underline v^B$ implies that \fer{log-xbar} holds indeed for all $(t,x)\in \R^+\times \R^d$ and hence the second inequality in \fer{uinf-lip} holds. The first inequality in \fer{uinf-lip} is a consequence of the second one.}\\

\noindent
(iii) We next prove that $\underline{v}^B(0,x) \geq u_0^B(x)=\max(u_0(x),-B)$, for all $x\in \R^d$. To this end, we first prove that 
\beq
\label{visco-ini}
\max( \p_t \underline v^B(0,x) - H(D_x \underline v^B(0,x)) -\widetilde R(x,\overline I(0)),\underline v^B(0,x)-u_0^B(x))\geq 0.
\eeq
To prove the above inequality, let $(\e_n,t_n,x_n)_n$ be such that, as $n\to +\infty$, $(\e_n,t_n,x_n)\to (0,0,x)$ and $v_{\e_n}^B(t_n,x_n)\to \underline v^B(0,x)$. Let's first suppose that there exists a subsequence, that we call again by an abuse of notation $(\e_n,t_n,x_n)_n$, such that $t_n=0$. It follows that 
$$
v_{\e_n}^B(t_n,x_n)=v_{\e_n}^{B}(0,x_n) = \e_n\log \big( e^{\f{u_{\e_n}^0(x_n)}{\e_n}}+e^{-\f{B}{\e}} \big).
$$
 We then let $n\to +\infty$ to obtain, thanks to \fer{as:ue0}, that $\underline v^B(0,x)=u_0^B(x)$ and hence \fer{visco-ini}.\\
 
\noindent
We now suppose that such {a} subsequence  does not exist and hence we can suppose that, removing if necessary a finite number of points from the sequence, for all $n\geq 1$, we have $t_n>0$. We can then repeat the arguments in Step (i) to prove that 
$$
 \p_t \underline v^B(0,x) - H(D_x \underline v^B(0,x)) -{\widetilde R(x,\overline I(0)) } \geq 0,
 $$
 and hence \fer{visco-ini}.\\
 
\noindent
We next prove that  $u_0^B ( x_0) \leq \underline v^B(0, x_0)$, following the arguments of \cite{GB:94}--Theorem 4.7. {To this end, we first notice that $H(p)\geq 0$, for all $p\in \R^d$,  since $e^{kp\cdot \nu}+e^{-kp\cdot \nu }\geq 2$.}\\

\noindent
 We consider the following test function
$$
\vp(t,x) = -\f{|x-x_0|^2}{\eta^2}-\f{t}{\eta}.
$$
For $\eta$ small enough, $\underline v^B -\vp$ attains a minimum at $(t_\eta,x_\eta)$ such that $t_\eta\to 0$ and $x_\eta\to x_0$ as $\eta \to 0$. Note that since $H$  
and $\widetilde R$ are bounded from below, for $\eta$ small enough,
$$
\p_t \vp(t_\eta,x_\eta) -H(D_x \vp(t_\eta,x_\eta))- {\widetilde R(x_\eta,\overline I(t_\eta))} <0.
$$
Since $\underline{v}^B$ is a supersolution to  \fer{eq:vB} for $t>0$, we deduce that $t_\eta=0$. Using  \fer{visco-ini}  we obtain that
$$
u_0^B(x_\eta) \leq \underline v^B(0,x_\eta).
$$
 Moreover, since $(0,x_\eta)$ is a minimum point of $\underline v^B -\vp$, we deduce that
 $$
u_0^B(x_\eta) \leq \underline v^B(0,x_\eta)\leq \underline v^B(0,x_0).
$$
Letting $\eta\to 0$, and thanks to the continuity of $u_0^B$ we obtain that
$$
u_0^B(x_0) \leq \underline v^B(0,x_0).
$$

(iv) We next prove that $\overline{v}^B(t,x)\leq \underline{v}^B(t,x)$. To this end, we first modify and regularize $\underline{v}^B(t,x)$ and then use the regularized function as a test function.\\
We first modify $\underline{v}^B$ at the initial time in the following way:
$$
\underline{v}^B_\diamond (t,x) =
\begin{cases}
\underline{v}^B (t,x) & t>0,\\
\displaystyle{\liminf_{\underset{s>0}{s\to 0}} \underline{v}^B(s,x)}& t=0.
\end{cases}
$$
Note from (iii) and the lower semi-continuity of   $\underline{v}^B$   that
$$
u_0^B(x)\leq \underline{v}^B_\diamond (0,x).
$$
Moreover, with this definition, {$\underline{v}^B_\diamond (t,x)$} is a viscosity supersolution of \fer{eq:vB} also on the boundary $t=0$.\\

\noindent
Note that thanks to \fer{uinf-lip} $\underline{v}^B_\diamond(t,x)$  is Lipschitz and a.e. differentiable with respect to $x$. We perform an inf-convolution of $\underline{v}^B_\diamond$ to make it also Lipschitz continuous with respect to time:
\beq
\label{infconv}
\underline{v}^B_{\diamond,\gamma} (t,x)= \inf_{s \in \R^+} \{ \underline{v}^B_\diamond(s,x)+ \f{(t-s)^2}{\gamma^2} \}.
\eeq
One can verify that $\underline{v}^B_{\diamond,\gamma}$ converges to $\underline{v}^B_{\diamond}$ as $\gamma \to 0$. We prove that $\underline{v}^B_{\diamond,\gamma}$ is a supersolution of a perturbed version of the equation in \fer{eq:vB} in $(0,+\infty)\times \R^d$. Let $\vp$ be a smooth test function and assume that $\underline{v}^B_{\diamond,\gamma}-\vp$ takes a minimum at the point $(t_0,x_0)\in (0,+\infty)\times \R^d$. Assume also that 
$$
\underline{v}^B_{\diamond,\gamma} (t_0,x_0)=\underline{v}^B_\diamond(s_0,x_0)+ \f{(t_0-s_0)^2}{\gamma^2}.
$$
Note that such $s_0\in [0,\infty)$ exists since the infimum in \fer{infconv} is attained.  Therefore, $(t_0,s_0,x_0)$ is a minimum point of the following function 
$$
(t,s,x)\mapsto   \underline{v}^B_\diamond(s,x)+ \f{(t-s)^2}{\gamma^2} -\vp(t,x).
$$
Since $\underline{v}^B_\diamond$ is a supersolution of \fer{eq:vB} in $ [0,\infty)\times \R^d$, we deduce that 
$$
\f{2(t_0-s_0)}{\gamma^2}-H( D_x  \vp(t_0,x_0))  - {\widetilde R(x_0, \overline I(s_0))} \geq 0,
$$
which is equivalent with 
$$
\p_t \vp (t_0,x_0)-H( D_x  \vp(t_0,x_0))  - {\widetilde R(x_0, \overline I(s_0))} \geq 0.
$$
We conclude that
\beq
\label{perturbed-sups}
\p_t \vp (t_0,x_0)-H( D_x  \vp(t_0,x_0))  - {\widetilde R(x_0,\overline I(t_0))} \geq  
{\widetilde R(x_0, \overline I(s_0(t_0)))}-  {\widetilde R(x_0, \overline I(t_0))},
\eeq
with $s_0(t_0)$ a   point where the infimum in \fer{infconv} is attained. Note also that, by the definition of $s_0(t_0)$ and the fact that $|\underline{v}^B_{\diamond}|$ is bounded,  there exists a constant $C$, {which may depend on $B$}, such that
\beq
\label{s0t0}
|t_0-s_0(t_0)| \leq C\gamma.
\eeq
However, despite the above inequality, the right hand side of \fer{perturbed-sups} can be large for small $\gamma$ because of the discontinuity of {$\overline I$}.
\\

 \noindent
Let $\chi_\beta(\cdot,\cdot)=\f{1}{\beta^{d+1}}\chi(\cdot/\beta,\cdot/\beta)$ be a   smoothing mollifier, with $\chi:\R^+\times \R^d\to \R^+$ a smooth function such that
$$
\begin{cases}
0\leq  \chi  \leq 1,\\
\int_{\R^+\times\R^d} \chi(t,x)dxdt =1,\\
\chi(t,x)=0, \qquad \text{if $|x|\geq 1$ or $|t|\geq 1$}.
\end{cases}
$$
We define
$$
\underline{v}_{\diamond ,\beta,\gamma}^B=\underline{v}^B_{\diamond,\gamma }  \ast \chi_\beta.
$$
Using the concavity of the Hamiltonian in \fer{eq:vB} {and \fer{perturbed-sups}} we obtain that
\beq
\label{HJ-reg}
\begin{array}{c}
\displaystyle \p_t \underline{v}_{\diamond,\beta,\gamma}^B (t_0,x_0)-\int_0^\infty \int_{\nu\in S^{d-1}}  \left( e^{k D_x  \underline{v}_{\diamond,\beta,\gamma}^B(t_0,x_0)\cdot \nu}-1 \right) \f{e^k dSdk}{|e^k-1|^{1+2\al}} - \widetilde R\ast \chi_\beta (t_0,x_0) \geq\\
\displaystyle\int_{\R^+\times \R^d} {\chi_\beta(t_0- \tau,x_0-y )} \big[ {\widetilde R(y, \overline I(s_0(\tau)))}-  {\widetilde R(y, \overline I(\tau))} \big] d\tau dy.
\end{array}
\eeq
We prove that the right hand side of the above inequality is greater than $-\f{\mu}{2(1-\mu)}$, with $0<\mu<1$ a small constant, for $\gamma$ small enough. To this end, define
$$
\mathcal A= \{ \tau \in [t_0-\beta, t_0 +\beta] \, | \,  | \, \overline I(s_0(\tau)) -  \overline I(\tau) \, | \leq \kappa \},
$$
with $\kappa$  a small constant to be chosen later. {Note that $\mathcal A$ may be empty.} We split the integral on the {r.h.s.} of \fer{HJ-reg} into two parts in the following way
$$
\begin{array}{c}
\displaystyle \int_ {\R^d} \int_{\R^+} {\chi_\beta(t_0- \tau,x_0-y )} \big[{ \widetilde R(y, \overline I(s_0(\tau)))}- { \widetilde R(y, \overline I(\tau))} \big] d\tau dy = \\
\displaystyle\int_ {y\in B_\beta(x_0)} \int_{\caa} {\chi_\beta(t_0- \tau,x_0-y )}  \big[{ \widetilde R(y, \overline I(s_0(\tau)))}- { \widetilde R(y, \overline I(\tau))} \big] d\tau dy+\\
\displaystyle\int_ {y\in B_\beta(x_0)} \int_{\caa^c}{\chi_\beta(t_0- \tau,x_0-y )}  \big[ {\widetilde R(y, \overline I(s_0(\tau)))}- { \widetilde R(y, \overline I(\tau))} \big] d\tau dy=
F_1+F_2.
\end{array}
$$
Using \fer{as:R2}, \fer{Gneg} and the definition of $\chi_\beta$ we obtain that
$$
| F_1 | \leq {K_1\kappa}.
$$
Moreover, using \fer{as:R3} and the definition of $\chi_\beta$ we obtain that, for some positive constant $C'$,
$$
| F_2 | \leq \f{C'K_2}{\beta}  \int_{\caa^c} d\tau.
$$
We then use the monotonicity of $\overline I$, {thanks to Proposition \ref{prop:reg},} \fer{boundI} and \fer{s0t0} to obtain that
$$
\int_{\caa^c} d\tau \leq \f{3C\gamma I_M}{\kappa},
$$
and hence 
$$
| F_2 | \leq \f{3C'CK_2 I_M\gamma}{\beta \kappa}   .
$$
Here, we have used the fact that at most $\f{I_M}{\kappa}$ disjoint intervals $[\tau_i,\tau_i+C\gamma]$ exist such that $I(\tau_i+C\gamma)-I(\tau_i) > \kappa$. Moreover, if $\mathcal B=\cup_{i\in \mathcal I} [\tau_i,\tau_i+C\gamma]$ is a maximal set of such intervals, then $\caa^c \subset \cup_{i\in \mathcal I} [\tau_i-C\gamma,\tau_i+2C\gamma]$.\\

\noindent
Combining the above properties and choosing $\kappa$ and $\gamma$ small enough such that $\gamma << \kappa \beta$, we obtain that
\beq
\label{HJ-reg-2}
\displaystyle \p_t \underline{v}_{\diamond,\beta,\gamma}^B (t_0,x_0)-\int_0^\infty \int_{\nu\in S^{d-1}}  \left( e^{k D_x  \underline{v}_{\diamond,\beta,\gamma}^B(t_0,x_0)\cdot \nu}-1 \right) \f{e^k dSdk}{|e^k-1|^{1+2\al}} - \widetilde R\ast \chi_\beta (t_0,x_0) \geq -\f{\mu}{2(1-\mu)}.
\eeq

\noindent
We thus obtain a supersolution, with a small error, which is smooth with respect to $x$ and  $t$. We then modify it to obtain a strict supersolution and also make the inequalities in \fer{uinf-lip} strict:
$$
\underline{v}_{\diamond,\beta,\gamma,\mu}^B =(1- \mu)   \underline{v}_{\diamond,\beta,\gamma}^B  + \mu  t, \quad \text{with $0<\mu<1$}.
$$
Using the concavity of the Hamiltonian  and the fact that $t$ is a strict supersolution of \fer{HJ-reg-2} we obtain that {$\underline{v}_{\diamond,\beta,\gamma,\mu}^B$} is a supersolution of the following equation, for $\gamma$ small enough,
\beq
\label{strict-super}
 \p_t \underline{v}_{\diamond, \beta,\gamma,\mu}^B -\int_0^\infty \int_{\nu\in S^{d-1}}  \left( e^{k D_x  \underline{v}_{\diamond,\beta,\gamma,\mu}^B\cdot \nu}-1 \right) \f{e^k dSdk}{|e^k-1|^{1+2\al}} \geq \widetilde R \ast \chi_\beta+\f\mu2, \quad \text{in $(0,+\infty)\times \R^d$}, 
\eeq
 and moreover 
\beq
\label{log-strict}
\|D_x \underline{v}_{\diamond ,\beta,\gamma,\mu}^B \|_{L^\infty(\R^d \times \R^+)}\leq 2\alpha (1-\mu),\qquad \underline{v}_{\diamond ,\beta,\gamma,\mu}^B(t,x+h)-\underline{v}_{\diamond, \beta,\mu}^B(t,x)\leq 2\alpha  (1-\mu)\log(1+|h|).
\eeq
Note also that by regularity and the above inequalities, $\underline{v}_{\diamond,\beta,\gamma,\mu}^B$ is a strict supersolution up to $t=0$.\\
We have now a good candidate for a test function. \\

\noindent
Fix $T>0$. Let's suppose that $\max_{(t,x)\in [0,T]\times \R^d}  \overline v^B - \underline v_{\diamond}^B\geq a >0$. Using the bound \fer{boundu} and the fact that $ \underline v^B\geq -B$, we obtain that such maximum is attained at some point $(t_0,x_0)\in [0,T]\times K$, with $K$ a compact set. Moreover, $\overline v^B(t_0,x_0)> -B$.  We can choose the set $K$ such that $x_0$ is an interior point of this set and such that the value of $\overline v^B - \underline v_{\diamond}^B$ on $[0,T]\times \p K$ is strictly less than this maximum.  
For $\gamma$, $\beta$ and $\mu$ small enough, $ \overline v^B -\underline{v}_{\diamond,\beta,\gamma,\mu}^B$  takes a positive maximum, greater than $a/2$, at some point $(\widetilde t   ,\widetilde x  )\in [0,T]\times K$, with $\overline v^B(\widetilde t,\widetilde x)> -B$. 
The main idea is to consider $\underline{v}_{\diamond,\beta,\gamma,\mu}^B$ as a test function at the point $(\widetilde t   ,\widetilde x  )$. To deal with the discontinuity in time of $\widetilde R$ we will use methods of viscosity solutions for Hamilton-Jacobi equations where the Hamiltonian is $L^1$ with respect to $t$ \cite{HI:85,PL.BP:87}. To this end, we define
$$
b_{\beta,\e}(t)= {\sup_{x\in K} \left({\widetilde R(x,I_\e(t))} - \big({\widetilde R(\cdot,I_\e(\cdot))} \ast \chi_\beta \big)(t,x) \right)},
$$ 
$$
b_{\beta}(t)={ \sup_{x \in K } \left( {\widetilde R(x,I(t))} - \big( {\widetilde R(\cdot,I(\cdot))} \ast \chi_\beta\big)(t,x) \right)}.
$$ 
One can verify that, for all $t>0$,
$$
\int_0^t b_{\beta,\e}(s)ds \to \int_0^t b_{\beta}(s)ds, \qquad \text{as $\e\to 0$},
$$
$$
 \int_0^t b_{\beta}(s)ds\to 0, \qquad \text{as $\beta\to 0$}.
 $$
 Therefore, for $\beta$ small enough, $\overline v^B(t,x) -\underline{v}_{\diamond,\beta,\gamma,\mu}^B (t,x) - \int_0^t b_{\beta}(s)ds$, attains a positive maximum.
Note from the definition of $\overline v^B$ and the above properties, there exists  a sequence $(\e_n)_n$, with $\e_n\to 0$ as $n\to \infty$,  such that $v_{\e_n}^B - \underline{v}_{\diamond,\beta,\gamma,\mu}^B-\int_0^t b_{\beta,\e_n}(s)ds$  takes a positive maximum at some point $(t_{\e_n},x_{\e_n} )\in K$. Passing to the limits along an appropriate subsequence, noting again by an abuse of notation $(\e_n)_{n}$, we obtain that, as $\e_n\to 0$, $(t_{\e_n} ,x_{\e_n})\to (\overline t  , \overline x   )$, such that $v_{\e_n}(t_{\e_n} ,x_{\e_n})\to \overline v^B(\overline t   ,\overline x   )$ and $(\overline t  , \overline x   )\in K$ is a maximum point of $\overline v^B(t,x) -\underline{v}_{\diamond,\beta,\gamma,\mu}^B (t,x) - \int_0^t b_{\beta}(s)ds$.\\

\noindent
Moreover, for all $(t,x) \in [0,T]\times K$, we have
\beq
\label{RbGe}
{\widetilde R(x,I_\e(t))} \leq b_{\beta,\e}(t)+{ \widetilde R(\cdot,I_\e(\cdot))} \ast \chi_\beta(t,x).
\eeq

\bigskip

{\bf Case $1:$} $\widetilde t >0$. Then, for $\e_n$   small enough, we have also  $ t_{\e_n }>0$.
 We then use $\overline v^B_{\diamond,\beta,\gamma,\mu}+\int_0^t b_{\beta,\e_n} (s)ds$ as a test function for equation \fer{eq:veB} on $v_{\e_n}$ at the point $(t_{\e_n },x_{\e_n} )$:
$$
 \begin{array}{rl}
\p_t  \underline{v}_{\diamond,\beta,\gamma,\mu}^B(t_{\e_n },x_{\e_n} ) &\leq \int_0^\infty \int_{\nu \in S^{d-1}} \Big(e^{\f{  \underline{v}_{\diamond,\beta,\gamma,\mu}^B(t_{\e_n},x_{\e_n}+(e^{\e_n k}-1)\nu) - \underline{v}_{\diamond,\beta,\gamma,\mu}^B(t_{\e_n },x_{\e_n })}{\e_n}} -1\Big)\f{e^k}{|e^k-1|^{1+2\alpha}}dS dk\\
& +\f{m_{\e_n}}{m_{\e_n}+e^{\f{-B}{\e_n}}}{ \widetilde R(x_{\e_n},I_{\e_n}(t_{\e_n}))}  - b_{\beta,\e_n}(t_{\e_n}) .
\end{array}
 $$
 Furthermore, thanks to \fer{RbGe} and taking $\e_n$, $\beta$, $\gamma$ and $\mu$ small enough such that $(t_{\e_n},x_{\e_n}) \in K$, we obtain
 $$
 \begin{array}{rl}
\p_t  \underline{v}_{\diamond,\beta,\gamma,\mu}^B(t_{\e_n },x_{\e_n} ) &\leq \int_0^\infty \int_{\nu \in S^{d-1}} \Big(e^{\f{  \underline{v}_{\diamond,\beta,\gamma,\mu}^B(t_{\e_n},x_{\e_n}+(e^{\e_n k}-1)\nu) - \underline{v}_{\diamond,\beta,\gamma,\mu}^B(t_{\e_n },x_{\e_n })}{\e_n}} -1\Big)\f{e^k}{|e^k-1|^{1+2\alpha}}dS dk\\
& +\f{m_{\e_n}}{m_{\e_n}+e^{\f{-B}{\e_n}}} {\widetilde R(x_{\e_n},I_{\e_n}(t_{\e_n}))}  - {\widetilde R(x_{\e_n},I_{\e_n}(t_{\e_n}))} +{\widetilde R(\cdot,I_{\e_n}(\cdot))} \ast \chi_\beta(t_{\e_n},x_{\e_n}).
\end{array}
 $$
 We then let $\e_n\to 0$ and use {similar arguments as in Step (i)}, \fer{log-strict} and the fact that $\overline v^B(\overline t,\overline x)> -B$  to find that 
\beq
\label{ineq-final}
 \begin{array}{c}
 \p_t \underline{v}_{\diamond, \beta,\gamma,\mu}^B(\overline t,\overline x) -\int_0^\infty \int_{\nu\in S^{d-1}}  \left( e^{k D_x  \underline{v}_{\diamond,\beta,\gamma,\mu}^B(\overline t, \overline x)\cdot \nu}-1 \right) \f{e^k dSdk}{|e^k-1|^{1+2\al}}\\
  \leq {\widetilde R(\cdot,I (\cdot))} \ast \chi_\beta(\overline  t,\overline x).
 \end{array}
\eeq
 {Note indeed that,  thanks to \fer{log-strict},
 $$
 \begin{array}{c}
 \int_M^\infty \int_{\nu \in S^{d-1}} \Big(e^{\f{  \underline{v}_{\diamond,\beta,\gamma,\mu}^B(t_{\e_n},x_{\e_n}+(e^{\e_n k}-1)\nu) - \underline{v}_{\diamond,\beta,\gamma,\mu}^B(t_{\e_n },x_{\e_n })}{\e_n}} -1\Big)\f{e^k}{|e^k-1|^{1+2\alpha}}dS dk\\
 \leq  \int_M^\infty \int_{\nu \in S^{d-1}} \Big(e^{2\al(1-\mu)k} -1\Big)\f{e^k}{|e^k-1|^{1+2\alpha}}dS dk
 \end{array}
 $$
which tends to $0$ as $M\to \infty$.}\\

\noindent 
{We conclude by noticing that} \fer{ineq-final} is in contradiction with  \fer{strict-super} which holds when $\gamma$ is chosen small enough.\\

{\bf Case $2:$} $\widetilde t=0$. If there is a subsequence  $(t_{\e_n } )_{\e_n }$ such that $t_{\e_n }=0$, then $\overline v^B(0,\widetilde x)=u_0^B(\widetilde x)$. This equality together with step (iii) implies that $ \overline v^B (0,\widetilde x)-\underline{v}_{\diamond,\beta,\gamma,\mu}^B (0,\widetilde x)< a/2$ for $\beta$, $\gamma$ and $\mu$ small enough,  which is a contradiction.
\\
We can thus assume that {$t_{\e_n}>0$}. Then, the problem can be treated as in Case $1$.

\section{Proof of Theorem \ref{th:main}}
\label{sec:proof-thm}

In view of Proposition \ref{prop:conv} which was proved in Section \ref{sec:proof-convergence}, to prove Theorem \ref{th:main} it remains to prove that $u$, a limit of $(u_\e)_\e$ along a subsequence as $\e\to 0$, satisfies \fer{constraint}, it is  minimal  in the set of viscosity supersolutions satisfying \fer{max} and it is a viscosity subsolution to \fer{HJ} in a weak sense as stated in the theorem. 
\\

\noindent {\bf [The proof of \fer{constraint}]} In view of Lemma \ref{lem:maxubar}, to prove \fer{constraint} it is enough to prove that $\max_{x\in \R^d} u(t,x)\leq 0$. This is immediate from the Hopf-Cole transformation  \fer{Hopf}, the Lipschitz continuity of $u$  and \fer{boundI}.\\

\noindent {\bf [Minimality of $u$ in the set of viscosity supersolutions satisfying \fer{max}]}   Let's assume that $\widetilde u$ is a viscosity supersolution to \fer{HJ} satisfying \fer{max}. To prove that $u \leq\widetilde u $ we first define analogously to Section \ref{sec:proof-convergence}:
$$
\widetilde v^B = \max (-B,\widetilde u-K_2t).
$$
We then repeat the arguments in the proof of Proposition \ref{lem:conv-vB}--Step (iv), to obtain that $v^B \leq \widetilde v^B$. We next let $B\to \infty$ to deduce that $u \leq\widetilde u$.
\\

\noindent {\bf [$u$ is a viscosity subsolution of \fer{HJ} in a weak sense]}
Let's suppose that $\vp\in \mathrm{C}( \R^+\times \R^d \to \R) \cap \mathrm{C}^2(\mathcal O(t_0,x_0))$ is a test function at the point $(t_0,x_0)$ such that it satisfies \fer{w-test-function} and such that $u -\vp$ takes a global and strict maximum at the point $(t_0,x_0)$. We prove that $\vp$ satisfies \fer{w-sub-cr}.\\

\noindent
Since $(u_\e)_\e$ converges locally uniformly to $u$, we deduce that $u_\e-\vp$ takes a global maximum at a point $(t_\e,x_\e) \in \R^+ \times \R^d$ such that, as $\e\to 0$, $(t_\e,x_\e) \to (t_0,x_0)$. Note also that replacing \fer{Hopf} in \fer{SMe} we obtain the following equation
$$
\p_t u_\e =\int_0^\infty \int_{\nu \in S^{d-1}} (e^{\f{u_\e (t,x+(e^{\e k}-1)\nu) - u_\e (t,x)}{\e}} -1) \f{e^k}{|e^k-1|^{1+2\alpha}}dS dk+  R(x,I_\e(t)).
$$
Thanks to the above equality and the fact that $u_\e-\vp$ takes a global maximum at the point $(t_\e,x_\e)$, we find that
$$
\p_t \vp (t_\e ,x_\e) \leq \int_0^\infty \int_{\nu \in S^{d-1}} (e^{\f{\vp (t_\e,x_\e+(e^{\e k}-1)\nu) - \vp (t_\e,x_\e)}{\e}} -1) \f{e^k}{|e^k-1|^{1+2\alpha}}dS dk+   R(x_\e,I_\e(t)).
$$
Note that for $\e$ small enough, $(t_\e,x_\e)\in B_r(t_0,x_0)$ with $r$ given by  \fer{w-test-function}. In view of  \fer{w-test-function}, we can pass to the limit in the above inequality and obtain  
$$
\p_t \vp(t_0,x_0) -\int_0^\infty \int_{\nu\in S^{d-1}}  \left( e^{k D_x  \vp(t_0,x_0)\cdot \nu}-1 \right) \f{e^k dSdk}{|e^k-1|^{1+2\al}} \leq \displaystyle\limsup_{ {\e\to 0}  } R(x_\e,I_\e(t_\e)),
$$
which leads to \fer{w-sub-cr} since $(t_\e,x_\e)\to (t_0,x_0)$ and thanks to the  estimate on $\f{d I_\e}{dt}$ (see Step (i) in the proof of Proposition \ref{lem:conv-vB} and \cite{GB.BP:08}). 
\\

\section{Proof of Theorem \ref{th:convn}}
\label{sec:proof-thm2}
Thanks to \fer{boundI} we obtain that, along subsequences as $\e\to 0$, $n_\e$ converges in {$\mathrm{L^\infty}\left( {w*}(0,\infty) ; \mathcal{M}^1(\R^d) \right)$} to a measure $n$. From the Hopf-Cole transformation and \fer{constraint}, \fer{supp-n} is immediate. We prove that \fer{zero-u-R} holds at the continuity points of $I(t)$.\\

\noindent 
Let $(t_0,x_0)$ be such that $u(t_0,x_0)=0$ with $t_0$ a continuity point of $I$. Then,  $\vp\equiv 0$ is a test function which satisfies \fer{w-test-function} and such that $u-\vp$ takes a maximum at the point $0$. Therefore, $\vp$ is an admissible test function and \fer{w-sub-cr} holds, i.e. 
$$
0 \leq \limsup_{s\to t_0} R(x_0,I(s))=R(x_0,I(t_0)).
$$
We next prove the inverse inequality. To this end, we integrate \fer{HJ} with respect to $t$, on $(t_0,t_0+h)$ at the point $x=x_0$. Using the positivity of
$$
\int_{S^{d-1}} \left( e^{k D_x  u\cdot \nu}-1 \right) dS \geq 0,
$$
and the fact that $u$ is a viscosity supersolution to \fer{HJ}, we obtain that
$$
 \f{1}{h} \int_{t_0}^{t_0+h} R(x_0,I(s))ds \leq \f{u(t_0+h,x_0)-u(t_0,x_0)}{h}.
$$
Using \fer{constraint} and the fact that $u(t_0,x_0)=0$ we obtain that
$$
 \f{1}{h} \int_{t_0}^{t_0+h} R(x_0,I(s))ds \leq 0.
$$
We then let $h\to 0$ and use the continuity of $I$ at the point $t_0$ to obtain
$$
R(x_0,I(t_0))\leq 0.
$$

\section{An example of a Hamilton-Jacobi equation of type \fer{HJ} with a solution not satisfying the second property of \fer{max}}
\label{sec:example}
In this section, we provide an example of a Hamilton-Jacobi equation of type \fer{HJ} which has a viscosity solution that does not satisfy the second property of \fer{max}. This example together with the fact that such Hamilton-Jacobi equation, with fixed $I$, has a unique viscosity solution (see \cite{CB.EC:13}, Section 6), indicates that the function $u$ might not be in general a viscosity solution of \fer{HJ}; it is only a viscosity supersolution and a viscosity subsolution in a weak sense as stated in Theorem \ref{th:main}.\\

\noindent
Consider the following equation 
\beq
\label{example}
\begin{cases}
\p_t u(t,x) -\int_0^\infty  \left( e^{k \p_x  u(t,x) } + e^{-k \p_x  u(t,x) } - 2 \right) \f{e^k dk}{|e^k-1|^{1+2\al}}  =a(t,x),\qquad(t,x)\in \R^+\times \R, \\
u(0,x)=0,\qquad x \in \R,
\end{cases}
\eeq
with 
$$
a(t,x) = \f{-C\sqrt{1+x^2}}{(1+t)^2} - \int_0^\infty \left(e^{\f{Ctx}{(1+t)\sqrt{1+x^2}}k}+e^{\f{-Ctx}{(1+t)\sqrt{1+x^2}}k}-2\right)  \f{e^k dk}{|e^k-1|^{1+2\al}},
$$
with $0<C< 2\alpha$. 
One can verify that 
$$
u(t,x)= -\f{ Ct\sqrt{1+x^2}}{1+t},
$$
is a solution to \fer{example}. However, $u$ does not satisfy the second property in  \fer{max} since it has linear decay.

\bigskip 

\textbf{Acknowledgments:} 
The author wishes to thank  Guy Barles and Emmanuel Chasseigne for very fruitful discussions.   The   author is also grateful for partial funding from the European Research Council (ERC)
under the European Union's Horizon 2020 research and innovation programme (grant agreement No
639638), held by Vincent Calvez, and from the french ANR projects KIBORD ANR-13-BS01-0004
and MODEVOL ANR-13-JS01-0009. 
 
%
\bibliographystyle{plain}

%

\end{document}